# INTERPRETING FORMULAS OF DIVISIBLE LATTICE ORDERED ABELIAN GROUPS - PRELIMINARY VERSION

MARCUS TRESSL

ABSTRACT. We show that a large class of divisible abelian $\ell$-groups (lattice ordered groups) of continuous functions is interpretable (in a certain sense) in the lattice of the zero sets of these functions. This has various applications to the model theory of these $\ell$-groups, including decidability results.

## Contents

## 1. Introduction

We show that a large class of divisible abelian $\ell$-groups of continuous functions is interpretable (in a certain sense, cf. 3.1) in the lattice of the zero sets of these functions. This interpretation will be independent of the respective group, and the interpretation will be recursive. The method (together with known results) has the following applications.

(1) If $K$ is an ordered field and $X \subseteq K^n$ is semilinear and of dimension 1, then the $\ell$-group $C_{s.l.}(X)$ of all $K$-semilinear continuous functions $X \longrightarrow K$ is decidable in the language $\{+, \leq\}$. This also implies decidability of the free vector lattice in 2 generators (as this is isomorphic to $C_{s.l.}(S^1)$ with base field $\mathbb{Q}$). This is a consequence of 4.2 (see also 4.3).
(2) The $\ell$-group of all continuous functions $\mathbb{R} \longrightarrow \mathbb{R}$ is decidable. This and the previous item rely on the decidability of the monadic second order theory of order, [Rabin1969]. See 4.4.
(3) The following $\ell$-groups are elementary equivalent (the language being $\{+, \leq\}$ again):
    (a) $C_{s.l.}(K)$ for any ordered field $K$

---

   (b) For any real closed field $R$, the $\ell$-group
   $$\{f : R \longrightarrow R \mid f \text{ is continuous and semi-algebraic}\}$$
   with addition.
   (c) For any real closed field $R$, the $\ell$-group
   $$\{f : R \longrightarrow (0, \infty) \subseteq R \mid f \text{ is continuous and semi-algebraic}\}$$
   with multiplication.

   This is done in 4.6 and 4.8 also for higher dimensions with the appropriate provisions.
   (4) If $M \prec N$ is an o-minimal extension of real closed fields and $n \in \mathbb{N}$, then the $\ell$-group of all continuous definable functions $M^n \longrightarrow M$ embeds elementarily (via the natural embedding) into the $\ell$-group of all continuous definable functions $N^n \longrightarrow N$. This relies on [Astier2013] and is done in 4.6.

After circulating an earlier version of this paper, I have been informed about the preliminary report [SheWei1987a] and the unpublished note [SheWei1987b], where similar applications to $\ell$-groups of continuous functions have been announced. An updated version of this paper is in preparation.

## 2. Constructible $\ell$-groups

**2.1. Definition.** Let $M$ be a divisible ordered abelian group (abbreviated DOAG), let $X$ be a Hausdorff space and let $C$ be a divisible $\ell$-group of continuous functions $X \longrightarrow M$. In other words, $C$ is a $\mathbb{Q}$-subvector lattice of the $\mathbb{Q}$-vector lattice $M^X$

(a) A subset $Z$ of $X$ is a $C$-**zero set** if $Z = \{f = 0\}$[1] for some $f \in C$. Complements of $C$-zero sets are called $C$-**cozero sets**. Boolean combinations of $C$-zero sets are called $C$-**constructible**.

(b) A $C$-**constructible function** is a map $g : Z \longrightarrow M$ for some $C$-constructible set $Z \subseteq X$ such that there are $C$-constructible sets $Z_1, \ldots, Z_n \subseteq Z$ with $Z = Z_1 \cup \ldots \cup Z_n$ and $f_1, \ldots, f_n \in C$ such that $g|_{Z_i} = f|_{Z_i}$ for all $i \in \{1, \ldots, n\}$.

We call $C$ a **constructible $\ell$-group** if the following conditions are satisfied:

(C1) (Zero set property)
Every closed and $C$-constructible set is a $C$-zero set.

(C2) (Closures)
The closure of a $C$-cozero set is $C$-constructible.

(C3) (Tietze extension)
If $A \subseteq X$ is closed and $C$-constructible, then every continuous and $C$-constructible function $A \longrightarrow M$ can be extended to a function from $C$.

(C4) (Glueing over closed sets)
If $A_1, \ldots, A_n \subseteq X$ are closed and $C$-constructible with $X = A_1 \cup \ldots \cup A_n$ and if $f : X \longrightarrow M$ is a function such that $f|_{A_i}$ is a continuous and $C$-constructible function for each $i$, then $f \in C$.

**2.2. *Examples.***
   (i) Let $M$ be an o-minimal expansion of a DOAG and let $X \subseteq M^n$ be $M$-definable. Let $C$ be the $\ell$-group of all $M$-definable and continuous functions $X \longrightarrow M$. Then properties (C2) and (C4) are clearly true for $C$.

---

[1] We write $\{f = 0\} = \{x \mid f(x) = 0\}$, $\{f < g\} = \{x \mid f(x) < g(x)\}$, etc.

Further, (C1) is true, because for every nonempty definable subset $A$ of $X$ that is closed in $X$, the distance function
$$\mathrm{dist}_A(\bar{x}) = \sup\{|x_i - a_i| \mid \bar{a} \in A,\ i \in \{1,\ldots,n\}\}$$
is in $C$ and has zero set $A$.

(ii) In the situation of (i), if $\dim(X) = 1$, and (C3) holds for 2-element subsets[2] then $C$ is a constructible $\ell$-group, since (C3) holds true using o-minimality. An interesting example here is when $M$ is an o-minimal expansion of the ordered $K$-vector space $K$ for some ordered field $K$.

(iii) Let $K$ be an ordered field and let $X \subseteq K^n$ be bounded and semilinear. Let $C$ be the $\ell$-group of all semilinear functions $X \longrightarrow K$. By [AscTha2015, Theorem 3.3], (C3) holds true in this case, hence $C$ is a constructible $\ell$-group (using (i)).

(iv) In the situation of (i), if $M$ expands a real closed field, then $C$ is a constructible $\ell$-group, since (C3) holds true in this case, [vdDries1998, Chapter 8, Cor. (3.10), p.138].

(v) Let again $M$ be an o-minimal expansion of a real closed field and let $X \subseteq M^n$ be $M$-definable. Let $C$ be the $\ell$-group of all $M$-definable functions $X \longrightarrow (0, +\infty)_M$ but this time the group law is pointwise multiplication. Then $C$ is a constructible $\ell$-group.

Notice that the $C$-zero sets are the sets of the form $\{f \neq 1\}$ with $f \in C$ and these are just the $M$-definable subsets of $X$ that are open in $X$. Properties (C1), (C2) and (C4) are true for the same reason stated in (i). Property (C3) holds true by the same reference to [vdDries1998, Chapter 8, Cor. (3.10), p.138].

(vi) If $X$ is a perfectly normal Hausdorff space (e.g. a metric space) and $C$ is the $\ell$-group of continuous functions $X \longrightarrow \mathbb{R}$, then $C$ is a constructible $\ell$-group.

Property (C3) holds true because the classical Tietze extension theorem holds in normal Hausdorff spaces. Property (C1) holds true by definition in a perfectly normal Hausdorff space. This also entails (C2). Property (C4) is obvious.

## 2.3. Construction of mixed functions to be used in 2.4

Let $X$ be a Hausdorff space and let $U \subseteq X$ be open. Let $A, B \subseteq X$ be closed with $A \cup B = X \setminus U$. Fix a divisible ordered abelian group $M$ (carrying the order topology) and suppose there is a continuous functions $d : X \longrightarrow M$ with zero set $A \cup B$.

Let $F, G, h : X \longrightarrow M$ be continuous such that

(a) $F \leq G$ and $F(u) < G(u)$ for all $u \in U$.
(b) $h(x) = F(x)$ for all $x \in A$ and $h(x) = G(x)$ for all $x \in B$. [3]

---

[2] By this we mean that for all $x \neq y \in X$ and all $a, b \in M$ there is some $f \in C$ with $f(x) = a$ and $f(y) = b$. An example where this fails is the following: Let $M \supsetneq (\mathbb{Q}, +, \leq)$ be a DOAG and assume $M$ has an element $b > \mathbb{Q}$. Take $X = [0, 1] \subseteq M$, $x = a = 0$ and $y = 1$. Then there is no $f \in C$ with $f(0) = 0$ and $f(1) = b$. The reason is that $f$ is piecewise linear, but the slope of each linear piece is rational; hence $f(1)$ is $\mathbb{Q}$-bounded

[3] When we apply this construction, we will just assume that $F(x) = G(x)$ for all $x \in A \cap B$ and then get $h$ from the Tietze extension property (C3) in the respective contexts, applied to the function $A \cup B \longrightarrow K$; $x \mapsto \begin{cases} F(x) & \text{if } x \in A, \\ G(x) & \text{if } x \in B. \end{cases}$

Then the function $H : X \longrightarrow M$ defined by

$$F + \frac{1}{2}\Big((|h - F| + |d|) \wedge (G - F)\Big) + \frac{1}{2}\Big((|h - F| + |d|) \wedge ((G - F - |d|) \vee 0)\Big)^4$$

extends $h|_{A \cup B}$ and satisfies $F(u) < H(u) < G(u)$ for all $u \in U$. (Notice that $H$ is contained in the divisible $\ell$-group generated by $F, G, h$ and $d$ in $M^X$.)

*Proof.* If $x \in A$, then $h(x) = F(x)$ and $d(x) = 0$, hence $(|h(x) - F(x)| + |d(x)|) = 0$ and so $H(x) = F(x) = h(x)$.

If $x \in B$, then $h(x) = G(x)$ and $d(x) = 0$, hence $(|h(x) - F(x)| + |d(x)|) = G(x) - F(x)$ and

$$H(x) = F(x) + \frac{1}{2}\Big((G(x) - F(x)) \wedge (G(x) - F(x))\Big)$$
$$+ \frac{1}{2}\Big((G(x) - F(x)) \wedge ((G(x) - F(x) - 0) \vee 0)\Big)$$
$$= F(x) + G(x) - F(x) = G(x) = h(x).$$

Now let $x \in U$. Then $|d(x)| > 0$ and as $F(x) < G(x)$ we have $(G(x) - F(x) - |d|(x)) \vee 0 < G(x) - F(x)$. Therefore

$$0 < \frac{1}{2}\Big((|h(x) - F(x)| + |d(x)|) \wedge (G(x) - F(x))\Big) \leq \frac{1}{2}(G(x) - F(x)),$$

and

$$0 \leq \frac{1}{2}\Big((|h(x) - F(x)| + |d(x)|) \wedge ((G(x) - F(x) - |d|(x)) \vee 0)\Big) < \frac{1}{2}(G(x) - F(x)).$$

Adding these up and adding $F(x)$ we get $F(x) < H(x) < G(x)$. □

We will now generalise (and correct) an argument of Amstislavskiy, see [Amstis2014, pp. 16-18] that he used to study the lattice of continuous functions $\mathbb{R} \longrightarrow \mathbb{R}$.

Given $n \in \mathbb{N}$ a finite family $(O_{ij})_{i,j \in I}$ of open subsets of a topological space $X$ is called an **open sign condition** if the following hold:

(O1) $O_{ii} = \emptyset$ for all $i \in I$.
(O2) $O_{ij} \cap O_{jk} \subseteq O_{ik}$ for all $i, j, k \in I$.
(O3) $O_{ik} \subseteq O_{ij} \cup O_{jk}$ for all $i, j, k \in I$. [5]

**2.4. Proposition.** *Let $M$ be a DOAG and let $C$ be a constructible $\ell$-group of continuous functions $X \longrightarrow M$ for some Hausdorff space $X$. Let $(O_{ij})_{i,j \in \{1,\ldots,n+k\}}$ be an open sign condition consisting of $C$-cozero sets and suppose there are $f_1, \ldots, f_n \in C$ such that*

$$O_{ij} = \{f_i < f_j\} \text{ for all } i, j \in \{1, \ldots, n\}.$$

*Then there are $f_{n+1}, \ldots, f_{n+k} \in C$ such that*

(∗) $\qquad\qquad O_{ij} = \{f_i < f_j\}$ *for all* $i, j \in \{1, \ldots, n + k\}.$

*Proof.* By a trivial induction on $k$ it suffices to do the case $k = 1$. First a couple of remarks and notations.

---

[4] Here and below $\wedge$ and $\vee$ denote infimum and supremum, respectively
[5] Modulo (O1) and (O2), it suffices to ask for $O_{ik} \subseteq O_{i,j} \cup O_{j,i} \cup O_{j,k} \cup O_{k,j}$

(a) For $x \in X$, let
$$I_x = \{i \in \{1,\ldots,n\} \mid x \in O_{i,n+1}\} \text{ and}$$
$$J_x = \{j \in \{1,\ldots,n\} \mid x \in O_{n+1,j}\}.$$
Then
  (1) $I_x \cap J_x = \emptyset$ otherwise for some $i \in \{1,\ldots,n\}$ we have $x \in O_{i,n+1} \cap O_{n+1,i} \subseteq O_{ii}$ by (ii), which contradicts (i).
  (2) For $x \in X$, $i \in I_x$ and $j \in J_x$ we have $x \in O_{i,n+1} \cap O_{n+1,j} \subseteq O_{ij}$, hence $f_i(x) < f_j(x)$.

(b) We define the closed and $C$-constructible set
$$E_i = \{x \in X \mid i \in \{1,\ldots,n\} \setminus (I_x \cup J_x)\} = X \setminus (O_{i,n+1} \cup O_{n+1,i}).$$
Further we define the closed and $C$-constructible set
$$E = E_1 \cup \ldots \cup E_n = X \setminus \bigcap_{i=1}^n (O_{i,n+1} \cup O_{n+1,i}).$$
By (O3), $f(x) := f_i(x)$ $(x \in E_i)$ gives a well defined function $f : E \longrightarrow M$ [6] and obviously $f$ is $C$-constructible and continuous.

(c) For any subset $I \subseteq \{1,\ldots,n\}$ we define
$$X_I := \{x \in X \mid I_x = I \text{ and } J_x = \{1,\ldots,n\} \setminus I\}.$$
Hence
$$X_I = \bigcap_{i \in I} O_{i,n+1} \cap \bigcap_{j \in \{1,\ldots,n\}\setminus I} O_{n+1,j}$$
is open and $C$-constructible and $X \setminus E = \bigcup_{I \subseteq \{1,\ldots,n\}} X_I$: the union is disjoint, since by (O2) and (O1) the intersection $O_{i,n+1} \cap O_{n+1,i}$ is empty for all $i$.

(d) For $\emptyset \neq I \subseteq \{1,\ldots,n\}$ we define $F_I : X \longrightarrow M$ by
$$F_I(x) = \max\{f_i(x) \mid i \in I\}.$$
For $I \subsetneq \{1,\ldots,n\}$ we define $H_I : X \longrightarrow M$ by
$$H_I(x) = \min\{f_i(x) \mid i \in \{1,\ldots,n\} \setminus I\}.$$
We extend this definition for the case $I = \emptyset$, or $I = \{1,\ldots,n\}$ as follows: Pick a function $c \in C$ that has no zeroes; such a function exists by (C1). Then define $F_\emptyset = H_\emptyset - |c|$ and $H_{\{1,\ldots,n\}} = F_{\{1,\ldots,n\}} + |c|$.

(e) The following hold true for all $I \subseteq \{1,\ldots,n\}$:
  (1) $F_I, H_I \in C$.
  (2) $F_I(x) < H_I(x)$ for all $x \in X_I$, by (a)(2).
  (3) Let $y \in X$ be a boundary point of $X_I$. Then $f(y) = F_I(y)$ or $f(y) = H_I(y)$.
      *Proof.* Since $y \in E = E_1 \cup \ldots \cup E_n$ there is some $i \in \{1,\ldots,n\}$ with $y \in E_i$ and we claim that
      - for $i \in I$ we have $f(y) = F_I(y)$.
      - for $i \in \{1,\ldots,n\} \setminus I$ we have $f(y) = H_I(y)$.

---

[6] If $i \in \{1,\ldots,n\} \setminus (I_x \cup J_x)$, then $(*)$ forces us to define $f_{n+1}(x) = f(x) = f_i(x)$.

To see this first notice that $f(y) = f_i(y)$, since $y \in E_i$.

First assume $i \in I$ and suppose $f_i(y) \neq F_I(y)$. Since $i \in I$ we have $f_i \leq F_I$, hence there is some $j \in I$ with $f_i(y) < f_j(y)$ and so $y \in O_{ij}$. Now (O3) entails $y \in O_{i,n+1} \cup O_{n+1,j}$. Since $y \in E_i$ we get $y \in O_{n+1,j}$. As $y$ is a boundary point of $X_I$ we know $X_I \cap O_{n+1,j} \neq \emptyset$. But this contradicts $j \in I$ and the definition of $X_I$.

Now assume $i \in \{1, \ldots, n\} \setminus I$ but $f_i(y) \neq H_I(y)$: Since $i \notin I$ we know $H_I \leq f_i$ and there is some $j \in \{1, \ldots, n\} \setminus I$ with $f_j(y) < f_i(y)$, hence $y \in O_{ji}$. Now (O3) entails $y \in O_{j,n+1} \cup O_{n+1,i}$. Since $y \in E_i$ we get $y \in O_{j,n+1}$. As $y$ is a boundary point of $X_I$ we get $X_I \cap O_{j,n+1} \neq \emptyset$. But this contradicts $j \notin I$. □

We are ready to define $f_{n+1}$. The aim is to extend $f$ (defined in (b)) to a function $f_{n+1} \in C$ such that $F_I(x) < f_{n+1}(x) < H_I(x)$ for all $x \in X_I$ and all $I \subseteq \{1, \ldots, n\}$ (This is actually forced by the assertion of the proposition).

By (C2) and (C4) it suffices to extend $f$ to a continuous $C$-constructible function $\overline{f} : \overline{X_I} \longrightarrow M$ with $F_I(x) < \overline{f}(x) < H_I(x)$ for all $x \in X_I$.

We apply 2.3 for the set $\overline{X_I}$ and the functions $F_I, H_I$ restricted to $\overline{X_I}$: We take $A = \{y \in \overline{X_I} \setminus X_I \mid f(y) = F_I(y)\}$ and $B = \{y \in \overline{X_I} \setminus X_I \mid f(y) = H_I(y)\}$. By (e)(3) we have $A \cup B = \overline{X_I} \setminus X_I$. By (C3) there is some $h \in C$ that extends $f$. By (C1) there is some $d \in C$ whose zero set is $A \cup B$. We have $F_I \leq H_I$ on $\overline{X_I}$ as is implied by (e)(2) and continuity of these functions. We see that all assumptions of 2.3 are satisfied and this yields the promised extension $\overline{f}$ of $f$ to $\overline{X_I}$.

This finishes the definition of $f_{n+1} \in C$. It remains to show that $O_{i,n+1} = \{f_i < f_{n+1}\}$ and $O_{n+1,i} = \{f_{n+1} < f_i\}$ for all $i$. For $x \in X$ we need to show that

(+) $f_i(x) < f_{n+1}(x) \iff x \in O_{i,n+1}$ and
(∗) $f_{n+1}(x) < f_i(x) \iff x \in O_{n+1,i}$.

If $x \notin E$ then $x \in X_I$ where $I = I_x$ and $F_I(x) < f_{n+1}(x) < H_I(x)$. By the definition of these functions we get (+) and (∗). Hence we may assume that $x \in E$. If $x \in E_i$, then $f_{n+1}(x) = f_i(x)$ and $x \notin O_{i,n+1} \cup O_{n+1,i}$ as required. Hence we may assume that $x \notin E_i$, thus

(†) $\qquad\qquad\qquad x \in O_{i,n+1} \cup O_{n+1,i},$

and $x \in E_j$ for some $j \neq i$, thus $f_{n+1}(x) = f_j(x)$ and

(‡) $\qquad\qquad\qquad x \notin O_{j,n+1} \cup O_{n+1,j}.$

If $f_i(x) < f_{n+1}(x)$, then $x \in O_{ij}$ and so $x \notin O_{n+1,j}$ forces $x \notin O_{n+1,i}$; thus $x \in O_{i,n+1}$ by (†).

Similarly, if $f_{n+1}(x) < f_i(x)$, then $x \in O_{ji}$ and so $x \notin O_{j,n+1}$ forces $x \notin O_{i,n+1}$; thus $x \in O_{n+1,i}$ by (†).

Since $O_{i,n+1}$ and $O_{n+1,i}$ are disjoint it remains to verify that $f_{n+1}(x) \neq f_i(x)$, i.e. $f_j(x) \neq f_i(x)$: This is implied by (O3) using (†) and (‡). □

## 3. Translating into the lattice of cozero sets

We work with formal linear combinations $t = q_1 x_1 + \ldots + q_n x_n$ with $q_i \in \mathbb{Q}$ in the variables $x_1, \ldots, x_n$ and call these expressions **terms**. We define the height of $t$ at $i$ to be $\text{height}_i(t) = \max\{|k|, |n|\}$, when $q_i = \frac{k}{n}$ is written in lowest terms and the **height of** $t$ as $\text{height}(t) = \max_1^n \text{height}_i(t)$. Note that $\text{height}(s + t) \leq$ $2 \cdot \text{height}(s) \cdot \text{height}(t)$ and $\text{height}(q \cdot t) \leq \text{height}(q) \cdot \text{height}(t)$ when $s, t$ are terms and $q \in \mathbb{Q}$. (However, notice that $\text{height}(qx_1 + qx_2) = \text{height}(q)$.)

We write $\text{Ht}_{\bar{x}}(d)$ for the set of terms in $x_1, \ldots, x_n$ of height at most $d$. This is extended to the empty tuple by setting $\text{Ht}_\emptyset(d) = \emptyset$ for all $d$. Observe that $\text{Ht}_{\bar{x}}(d)$ is a finite set. For $f_1, \ldots, f_{n+k} \in C$ and $t = q_1 x_1 + \ldots + q_n x_n \in \text{Ht}_{\bar{x}}(d)$ we write $t(f_1, \ldots, f_{n+k}) = q_1 f_1 + \ldots + q_n f_n \in C$.

Let $M$ be a DOAG and let $C$ be a constructible $\ell$-group of continuous functions $X \longrightarrow M$ for some Hausdorff space $X$ (cf. 2.1). Let $L_C$ be the partially ordered set of $C$-cozero sets. Let $(f_1, \ldots, f_n) \in C^n$ and let $\bar{x}$ be an $n$-tuple of variables. We define
$$\Lambda_{d,\bar{x}}(\bar{f}) = (\{s(\bar{f}) < t(\bar{f})\})_{s,t \in \text{Ht}_{\bar{x}}(d)},$$
which is a finite tuple with entries in $L_C$. Notice that $\Lambda_{d,\bar{x}}(\bar{f})$ depends on $\bar{x}$ only in the choice of the index set enumerating the tuple $(\{s(\bar{f}) < t(\bar{f})\})_{s,t \in \text{Ht}_{\bar{x}}(d)}$. However, later on we will need to know this index set.

3.1. We intend to define for each $\{\leq, +, -, 0\}$-formula $\varphi(x_1, \ldots, x_n)$ a $\{\leq\}$-formula $\varphi_\Lambda$ together with a number $d \in \mathbb{N}$ (depending on $\varphi$), such that the number of free variables of $\varphi_\Lambda$ is at most $\left(\#\text{Ht}_{\bar{x}}(d)\right)^2$ and such that for every constructible $\ell$-group of functions $X \longrightarrow M$ and all $\bar{f} \in C^n$ we have
$$C \models \varphi(\bar{f}) \iff L_C \models \varphi_\Lambda(\Lambda_{d,\bar{x}}(\bar{f})).$$
In this sense the formula $\varphi$ is interpreted in $L_C$, explaining the title of the article. However, the method does not interpret $C$ in $L_C$ in the standard model theoretic sense. For example, $X$ could be a singleton. Then $L_C$ is the 2-element lattice and $C$ is a DOAG contained in $M$.

The definition of $\varphi_\Lambda$ will be by induction on the complexity of $\varphi$ and the case of an existential quantifier is the central step. We briefly explain the idea. Suppose $\varphi(\bar{x}, y)$, $\bar{x} = (x_1, \ldots, x_n)$ is a quantifier free $\{\leq, +, -, 0\}$-formula and $\bar{f}$ is an $n$-tuple from some constructible $\ell$-group $C$. We want to write down a condition on $C$-cozero sets determined by $\bar{f}$, which guarantees the existence of some $g \in C$ with $C \models \varphi(\bar{f}, g)$. Suppose for a moment, $\varphi$ would just be a formula in the language $\{\leq\}$. Then 2.4 gives an answer by saying:

> Find a sign condition $(O_{ij})_{i,j \in \{1,\ldots,n+1\}}$ that extends the sign condition $(\{f_i < f_j\})_{i,j \in \{1,\ldots,n\}}$ given by $\bar{f}$.

(Notice that $\varphi$ is currently just a Boolean combination of inequalities $x_i \leq x_j$ and $y \leq x_i$, etc.; hence if we evaluate $\varphi$ we get $C \models g \leq f_i \iff \{f_i < g\} = \emptyset$. Thus, realizing the sign condition $(O_{ij})_{i,j \in \{1,\ldots,n+1\}}$ as $(\{f_i < f_j\})_{i,j \in \{1,\ldots,n+1\}}$ with $f_{n+1} := g$ is the same as finding a solution $g$ of $\varphi(\bar{f}, g)$ in $C$.)

In general, the idea now is to consider the quantifier free $\{\leq, +, -, 0\}$-formula $\varphi$ as a formula in the language $\{\leq\}$ where the variables are replaced by terms $t(\bar{x}, y)$. The terms $t(\bar{x}, y)$ are of the form $ky + s(\bar{x})$ with terms $s(\bar{x})$ in $\bar{x}$ and $k \in \mathbb{Z}$ (we will work with $k \in \mathbb{Q}$ right away). Since we are working with divisible groups we may then also assume that all occurrences of the terms $t(\bar{x}, y)$ in atomic subformulas of $\varphi$ separate the variables, i.e. are of the form $y \leq s(\bar{x})$, or $y \geq s(\bar{x})$. We can then apply the idea above to the long sign condition given by the $s(\bar{f})$ and get the required solution $g$. This central step is carried out in 3.4. The formulation of the condition in terms of first order logic, addressing an inductive assumption on all formulas $\varphi$ (not only quantifier free ones) is done in 3.5. The ingredients are then assembled in 3.8.

In order to carry this out it is convenient to work with the following set-up.

3.2. *Convention.* In what follows, $\mathscr{L}^+$ stands for the first order language

$$\{+, -, (f_q \mid q \in \mathbb{Q}), \leq\},$$

where $+, -$ are binary function symbols, $f_q$ are unary function symbols and $\leq$ is a binary relation symbol. The variables in this language are denoted by $v_n$ with $n \in \mathbb{N}$. If we talk about tuples of variables $\bar{x}$ below, then this always means a finite tuple $(v_{i_1}, \ldots, v_{i_n})$, where $i_j \neq i_k$ for $j \neq k$. The $\mathscr{L}^+$-structures that we are interested in are divisible $\ell$-groups and the symbols will have their natural interpretation ($f_q(x)$ is scalar multiplication of $x$ by $q$, which will of course be written as $q \cdot x$). Thus, a term as described at the beginning of this section is just a term in $\mathscr{L}^+$.

Further, the language $\mathscr{L}^{lat}$ (of posets with bottom element) is the first order language $\{\leq, \bot\}$, where the variables are denoted by $v_{s,t}$ with $\mathscr{L}^+$-terms $s, t$. Tuples of distinct variables in this language will be abbreviated by capital letters $\bar{X}, \bar{Y}$, etc.

In $\mathscr{L}^+$ we will also use the abbreviation $a \wedge b$ and $a \vee b$ for the infimum and supremum of $\{a, b\}$ for $\leq$.

Next we define auxiliary formulas used in the definition of the $\varphi_\Lambda$, see 3.3. To explain these, let $M$ be a DOAG and consider a constructible $\ell$-group $C$ of continuous functions $X \longrightarrow M$ for some Hausdorff space $X$ (cf. 2.1). Let $(f_1, \ldots, f_n) \in C^n$ and let $\bar{x}$ be an $n$-tuple of variables. For any $d \in \mathbb{N}$ the family $(O_{s,t})_{s,t \in \mathrm{Ht}_{\bar{x}}(d)}$ of open subsets of $X$ defined by

$$O_{s,t} = \{s(f_1, \ldots, f_n) < t(f_1, \ldots, f_n)\}$$

satisfies (O1), (O2) and (O3) (for the index set $\mathrm{Ht}_{\bar{x}}(d)$) as well as

(O4) $O_{s+r,t+r} = O_{s,t}$ for all $r, s, t \in \mathrm{Ht}_{\bar{x}}(d)$ with $s+r, t+r \in \mathrm{Ht}_{\bar{x}}(d)$.
(O5) $O_{r+s,t} = O_{r,t-s}$ for all $r, s, t \in \mathrm{Ht}_{\bar{x}}(d)$ with $r+s, t-s \in \mathrm{Ht}_{\bar{x}}(d)$.
(O6) $O_{q \cdot s, t} = O_{s, \frac{1}{q} \cdot t}$ for all $s, t \in \mathrm{Ht}_{\bar{x}}(d)$, $q \in \mathbb{Q}$, $q > 0$ with $s, t, q \cdot s, \frac{1}{q} \cdot t \in \mathrm{Ht}_{\bar{x}}(d)$.

3.3. **Definition of $\delta_{d,\bar{x}}$, $\Delta_{d,\bar{x},y}$ and $D_\varphi$**

(a) For $d \in \mathbb{N}$ and a tuple $\bar{x} = (x_1, \ldots, x_n)$ of variables let $\delta_{d,\bar{x}}$ be the following formula in the language $\mathscr{L}^{lat}$:

$\delta_{d,\bar{x}}$ is the $\mathscr{L}^{lat}$-formula in the variables $v_{s,t}$ with $s, t \in \mathrm{Ht}_{\bar{x}}(d)$ expressing that (O1)-(O6) hold for $\mathrm{Ht}_{\bar{x}}(d)$. Explicitly, $\delta_{d,\bar{x}}$ is the conjunction of the following $\mathscr{L}^{lat}$-formulas:
  (o1) $v_{s,s} = \bot$ for all $s \in \mathrm{Ht}_{\bar{x}}(d)$
  (o2) $v_{r,s} \wedge v_{s,t} \leq v_{r,t}$ for all $r, s, t \in \mathrm{Ht}_{\bar{x}}(d)$
  (o3) $v_{r,t} \leq v_{r,s} \vee v_{s,t}$ for all $r, s, t \in \mathrm{Ht}_{\bar{x}}(d)$
  (o4) $v_{s+r,t+r} = v_{s,t}$ for all $r, s, t \in \mathrm{Ht}_{\bar{x}}(d)$ with $s+r, t+r \in \mathrm{Ht}_{\bar{x}}(d)$.
  (o5) $v_{r+s,t} = v_{r,t-s}$ for all $r, s, t \in \mathrm{Ht}_{\bar{x}}(d)$ with $r+s, t-s \in \mathrm{Ht}_{\bar{x}}(d)$.
  (o6) $v_{q \cdot s, t} = v_{s, \frac{1}{q} \cdot t}$ for all $s, t \in \mathrm{Ht}_{\bar{x}}(d)$, $q \in \mathbb{Q}$, $q > 0$ with $s, t, q \cdot s, \frac{1}{q} \cdot t \in \mathrm{Ht}_{\bar{x}}(d)$.

(b) For $d \in \mathbb{N}$ and a tuple $\bar{x} = (x_1, \ldots, x_n)$ of variables and another single variable $y$ let $\Delta_{d,\bar{x},y}$ be the following formula in the language $\mathscr{L}^{lat}$. Let $\bar{W}$ be the tuple of all the variables $v_{s,t}$ with $s,t \in \mathrm{Ht}_{(\bar{x},y)}((2d)^4)$ such that

$$\{s,t\} \not\subseteq \mathrm{Ht}_{\bar{x}}((2d)^4) \text{ and } \{s,t\} \not\subseteq \mathrm{Ht}_{(\bar{x},y)}(d).$$

Then $\Delta_{d,\bar{x},y}$ is the $\mathscr{L}^{lat}$-formula $\exists \bar{W}\ \delta_{(2d)^4,(\bar{x},y)}$. Hence the free variables of $\Delta_{d,\bar{x},y}$ are exactly those $v_{s,t}$ with $\{s,t\} \subseteq \mathrm{Ht}_{\bar{x}}((2d)^4)$, or, $\{s,t\} \subseteq \mathrm{Ht}_{(\bar{x},y)}(d)$. Further, $\Delta_{d,\bar{x},y}$ is an existential formula, since $\delta$ is quantifier free.

(c) For an $\mathscr{L}^+$-formula $\varphi$, let $d_\varphi$ be the maximum of all heights of terms appearing in $\varphi$, let $N_\varphi$ be the number of quantifiers appearing in $\varphi$ and let $D_\varphi$ be the result of $N_\varphi$ iterations of the function $(2x)^4$ starting with $d_\varphi$.[7]

For the inductive definition of the $\varphi_\Lambda$ we first need two lemmas, which will allow us to process existential quantifiers.

**3.4. Lemma.** *Let $M$ be a DOAG and consider a constructible $\ell$-group $C$ of continuous functions $X \longrightarrow M$ for some Hausdorff space $X$. Let $\bar{x}$ be an $n$-tuple of variables and let $y$ be another single variable. Suppose $d \in \mathbb{N}$, $f_1, \ldots, f_n \in C$ and we are given a family*

$$(O_{s,t})_{s,t \in \mathrm{Ht}_{(\bar{x},y)}((2d)^4)}$$

*of $C$-cozero sets of $X$ satisfying (O1)–(O6) for $\mathrm{Ht}_{(\bar{x},y)}((2d)^4)$ such that*

(†) $$O_{s,t} = \{s(\bar{f}) < t(\bar{f})\} \text{ for all } s,t \in \mathrm{Ht}_{\bar{x}}((2d)^4).$$

*Then there is some $g \in C$ with*

$$O_{s,t} = \{s(\bar{f},g) < t(\bar{f},g)\} \text{ for all } s,t \in \mathrm{Ht}_{(\bar{x},y)}(d).$$

*Proof.* By 2.4 we know that there are $g_t \in C$ for all $t \in \mathrm{Ht}_{(\bar{x},y)}((2d)^4)$ such that $O_{s,t} = \{g_s < g_t\}$ for all $s,t \in \mathrm{Ht}_{(\bar{x},y)}((2d)^4)$ and such that

(‡) $$g_s = s(\bar{f}) \text{ for all } s \in \mathrm{Ht}_{\bar{x}}((2d)^4).$$

We claim that $g := g_y$ has the required properties.

Take $s,t \in \mathrm{Ht}_{(\bar{x},y)}(d)$. Hence there are $s_0, t_0 \in \mathrm{Ht}_{\bar{x}}(d)$ and $p,q \in \mathbb{Q}$ of height at most $d$ such that $s = qy + s_0$ and $t = py + t_0$. Then $q - p$ is of height at most $2d^2 \leq (2d)^2$ and so $(q-p)y + s_0 \in \mathrm{Ht}_{(\bar{x},y)}((2d)^2)$. If $p = q$ then the assertion follows

---

[7]Explicitly, $D_\varphi = 2^{\frac{4^{N+1}-4}{3}} \cdot d_\varphi^{4^N}$ with $N = N_\varphi$

easily from (O4) and the properties of the $g_r$. Suppose $q > p$. Then

$$\begin{aligned}
O_{s,t} &= O_{qy+s_0, py+t_0} \\
&= O_{(q-p)y+s_0, t_0} \text{ by (O4) since } (q-p)y + s_0 \in \text{Ht}_{(\bar{x},y)}((2d)^2) \\
&= O_{(q-p)y, t_0-s_0} \text{ by (O5) since } (q-p)y, t_0 - s_0 \in \text{Ht}_{(\bar{x},y)}((2d)^2) \\
&= O_{y, \frac{1}{q-p}(t_0-s_0)} \text{ by (O6) since } q > p \text{ and } \frac{1}{q-p}(t_0 - s_0) \in \text{Ht}_{(\bar{x},y)}((2d)^4) \\
&= \{g_y < g_{\frac{1}{q-p}(t_0-s_0)}\} \text{ by choice of the } g_r \\
&= \{g_y < \frac{1}{q-p}(t_0-s_0)(\bar{f})\} \text{ by ($\ddagger$) and as } \frac{1}{q-p}(t_0-s_0) \in \text{Ht}_{\bar{x}}((2d)^4) \\
&= \{g < \frac{1}{q-p}(t_0(\bar{f}) - s_0(\bar{f}))\} \text{ by definition of } g \\
&= \{(q-p)g < t_0(\bar{f}) - s_0(\bar{f})\} \text{ since } q > p \\
&= \{qg + s_0(\bar{f}) < pg + t_0(\bar{f})\} \\
&= \{s(\bar{f}, g) < t(\bar{f}, g)\}.
\end{aligned}$$

If $q < p$ we can do the same reasoning but this time we start with $O_{s,t} = O_{s_0, (p-q)y+t_0}$. □

**3.5. Lemma.** *Take an $\mathscr{L}^+$-formula $\varphi(\bar{x}, y)$, where $\bar{x}$ is an $n$-tuple and $y$ is a single variable. Let $d \geq d_\varphi$. Suppose $\varphi_\Lambda$ is an $\mathscr{L}^{lat}$-formula in at most the free variables $v_{s,t}$, where $s,t \in \text{Ht}_{(\bar{x},y)}(d)$. We define an $\mathscr{L}^{lat}$-formula $\Phi$[8] as*

$$\exists \bar{Y} \left( \Delta_{d,\bar{x},y}(\bar{X}, \bar{U}, \bar{Y}) \ \& \ \varphi_\Lambda(\bar{X}, \bar{Y}) \right),$$

*where $\bar{X}$ is the tuple of all variables $v_{s,t}$ with $s,t \in \text{Ht}_{\bar{x}}(d)$, $\bar{Y}$ is the tuple of the variables $v_{s,t}$ with $s,t \in \text{Ht}_{(\bar{x},y)}(d)$ such that $\{s,t\} \not\subseteq \text{Ht}_{\bar{x}}(d)$ and $\bar{U}$ is the tuple of all variables $v_{s,t}$ with $s,t \in \text{Ht}_{\bar{x}}((2d)^4)$ such that $\{s,t\} \not\subseteq \text{Ht}_{\bar{x}}(d)$.*

*If $C$ is a constructible $\ell$-group of functions $X \longrightarrow M$ such that for all $(\bar{f}, g) \in C^{n+1}$ we have*

$$C \models \varphi(\bar{f}, g) \iff L_C \models \varphi_\Lambda(\Lambda_{d,(\bar{x},y)}(\bar{f}, g)),$$

*then for all $\bar{f} \in C^n$ we have*

$$C \models (\exists y \, \varphi)(\bar{f}) \iff L_C \models \Phi(\Lambda_{(2d)^4, \bar{x}}(\bar{f})).$$

*Proof.*

$$\begin{aligned}
C \models (\exists y \, \varphi)(\bar{f}) &\iff \text{there is } g \in C \text{ with } C \models \varphi(\bar{f}, g) \\
&\iff \text{there is } g \in C \text{ with } L_C \models \varphi_\Lambda(\Lambda_{d,(\bar{x},y)}(\bar{f}, g)), \text{ by assumption} \\
&\iff \text{there is a } \bar{Y}\text{-tuple } \bar{O} \subseteq L \text{ with} \\
& \quad L_C \models \Delta_{d,\bar{x},y}(\Lambda_{(2d)^4, \bar{x}}(\bar{f}), \bar{O}) \ \& \ \varphi_\Lambda(\Lambda_{d,\bar{x}}(\bar{f}), \bar{O}),
\end{aligned}$$

where $\bar{Y}$ is the tuple of the variables $v_{s,t}$ with $s,t \in \text{Ht}_{(\bar{x},y)}(d)$ and $\{s,t\} \not\subseteq \text{Ht}_{\bar{x}}(d)$.

---

[8] One should think of $\Phi$ as $(\exists y \, \varphi)_\Lambda$ once $d$ is chosen appropriately. Notice that $\Phi$ is a sentence if $\varphi$ is a formula in the variable $y$

*Proof of the last equivalence.*

$\Rightarrow$. Take $g \in C$ with $L_C \models \varphi_\Lambda(\Lambda_{d,(\bar{x},y)}(\bar{f},g))$. For $s,t \in \mathrm{Ht}_{(\bar{x},y)}((2d)^4)$ let $O_{s,t} = \{s(\bar{f},g) < t(\bar{f},g)\}$. The resulting family $(O_{s,t})_{s,t \in \mathrm{Ht}_{(\bar{x},y)}((2d)^4)}$ satisfies (O1)-(O6), hence $L_C \models \delta_{(2d)^4,(\bar{x},y)}((O_{s,t})_{s,t \in \mathrm{Ht}_{(\bar{x},y)}((2d)^4)})$. Hence if $\bar{O}$ is the enumeration of all the $O_{s,t}$ with $s,t \in \mathrm{Ht}_{(\bar{x},y)}(d)$ such that $\{s,t\} \not\subseteq \mathrm{Ht}_{\bar{x}}(d)$, then

$$L_C \models \Delta_{d,\bar{x},y}(\Lambda_{(2d)^4,\bar{x}}(\bar{f}),\bar{O}) \,\&\, \varphi_\Lambda(\Lambda_{d,\bar{x}}(\bar{f}),\bar{O})$$

(notice that $\Lambda_{d,\bar{x},y}(\bar{f},g) = (\Lambda_{d,\bar{x}}(\bar{f}),\bar{O})$).

$\Leftarrow$. Take a $\bar{Y}$-tuple $\bar{O} \subseteq L_C$ with $L_C \models \Delta_{d,\bar{x},y}(\Lambda_{(2d)^4,\bar{x}}(\bar{f}),\bar{O}) \,\&\, \varphi_\Lambda(\Lambda_{d,\bar{x}}(\bar{f}),\bar{O})$. Hence there is a family $(O_{s,t})_{s,t \in \mathrm{Ht}_{(\bar{x},y)}((2d)^4)}$ of elements of $L_C$ satisfying (O1)-(O6) such that $O_{s,t} = \{s(\bar{f}) < t(\bar{f})\}$ for all $s,t \in \mathrm{Ht}_{\bar{x}}((2d)^4)$ and such that $\bar{O}$ is the enumeration of the $O_{s,t}$ with $s,t \in \mathrm{Ht}_{(\bar{x},y)}(d)$ such that $\{s,t\} \not\subseteq \mathrm{Ht}_{\bar{x}}(d)$. By 3.4 there is some $g \in C$ with

$$O_{s,t} = \{s(\bar{f},g) < t(\bar{f},g)\} \text{ for all } s,t \in \mathrm{Ht}_{(\bar{x},y)}(d).$$

But this means $L_C \models \varphi_\Lambda(\Lambda_{d,(\bar{x},y)}(\bar{f},g))$, as required. $\square$

Hence the equivalence is justified and we see that $\Phi$ has the required property. $\square$

**3.6. Definition.** Now for the definition of the $\mathscr{L}^{lat}$-formula $\varphi_\Lambda$. Let $\varphi$ be an $\mathscr{L}^+$-formula.

(I) If $\varphi$ is $s \leq t$ with $\mathscr{L}^+$-terms $s,t$, then $\varphi_\Lambda$ is $v_{t,s} \doteq \bot$.
(II) $(\varphi \,\&\, \psi)_\Lambda$ is $\varphi_\Lambda \,\&\, \psi_\Lambda$.
(III) $(\neg \varphi)_\Lambda$ is $\neg(\varphi_\Lambda)$.
(IV) Let $\varphi$ be an $\mathscr{L}^+$-formula with exactly the free variables $v_{i_1},\ldots,v_{i_m}$, where $i_1 < \ldots < i_m$ and let $y$ be any variable. If $y$ is not among the $v_{i_j}$, then we define $(\exists y\, \varphi)_\Lambda$ as $\varphi_\Lambda$. If $y = v_{i_j}$ for some $j$, then we define $(\exists y\, \varphi)_\Lambda$ as the formula $\Phi$ from 3.5 for the choice $\bar{x} = (v_{i_1},\ldots,v_{i_{j-1}},v_{i_{j+1}},\ldots,v_{i_m})$ and $d = D_\varphi$ (defined in 3.3(c)).

**3.7. Lemma.** *If the free variables of $\varphi$ are among $\bar{x}$, then the free variables of $\varphi_\Lambda$ are among the variables $v_{s,t}$, where $s,t \in \mathrm{Ht}_{\bar{x}}(D_\varphi)$.*

*Proof.* This is obvious if $\varphi$ is $s \leq t$ and the induction in the cases (II) and (III) is clear, too.

For case (IV) we refer to 3.5. Let $\varphi$ be an $\mathscr{L}^+$-formula with exactly the free variables $v_{i_1},\ldots,v_{i_m}$, where $i_1 < \ldots < i_m$ and let $y$ be any variable. If $y$ is not among the $v_{i_j}$, then $(\exists y\, \varphi)_\Lambda$ is $\varphi_\Lambda$ and we are done.

So assume $y = v_{i_j}$ for some $j$ and take $\Phi$ from 3.5 for the choice $d = D_\varphi \geq d_\varphi$ ($= d_{\exists y \varphi}$) and $\bar{x} = (v_{i_1},\ldots,v_{i_{j-1}},v_{i_{j+1}},\ldots,v_{i_m})$. By induction we know that the free variables of $\varphi_\Lambda$ are among the variables $v_{s,t}$, where $s,t \in \mathrm{Ht}_{(\bar{x},y)}(D_\varphi)$.

Then by 3.5 the free variables of $\Phi$ are all $v_{s,t}$ with $s,t \in \mathrm{Ht}_{\bar{x}}(d)$ and all the variables $v_{s,t}$ with $s,t \in \mathrm{Ht}_{\bar{x}}((2d)^4)$ such that $\{s,t\} \not\subseteq \mathrm{Ht}_{\bar{x}}(d)$. Since $(2d)^4 = D_{\exists y\, \varphi}$, this shows the claim. $\square$

**3.8. Theorem.** *Let $C$ be a constructible $\ell$-group of functions $X \longrightarrow M$ and let $L$ be the poset of $C$-cozero sets of $C$. Let $\varphi$ be an $\mathscr{L}^+$-formula in at most the free variables $\bar{x}$. Let $d = D_\varphi$ as defined in 3.3(c). Then for all $\bar{x}$-tuples $\bar{f}$ with entries in $C$ we have*

$$C \models \varphi(\bar{f}) \iff L_C \models \varphi_\Lambda(\Lambda_{d,\bar{x}}(\bar{f})).^9$$

*Further, if $\varphi$ is an $\forall_n$-formula, then so is $\varphi_\Lambda$, and similarly for $\exists_n$.*

*Proof.* By induction on the complexity of $\varphi$. To start with, let $\varphi$ be $s_0 \leq t_0$, where $s_0, t_0$ are terms whose variables are among $x_1, \ldots, x_n$. We need to show that

$$s_0(\bar{f}) \leq t_0(\bar{f}) \iff L_C \models v_{t_0,s_0} \doteq \bot \ (\Lambda_{d,\bar{x}}(\bar{f})),$$

where $d = d_\varphi$ (as $q = 0$ for this formula). Now plugging $\Lambda_{d,\bar{x}}(\bar{f})$ into the formula $v_{t_0,s_0} \doteq \bot$ is the same as plugging $\{t_0(\bar{f}) < s_0(\bar{f})\} \in L_C$ into that formula. Hence the equivalence holds true.

In the cases (II) and (III) the induction is clear. Case (IV) is precisely covered by 3.5.

The statement about the quantifier complexity follows by induction on this complexity from the definition of $\varphi_\Lambda$ using the fact that $\Delta_{d,\bar{x},y}(\bar{X}, \bar{U}, \bar{Y})$ in 3.5 is an existential formula (cf. 3.3(b)). $\square$

**3.9. Corollary.** *In the situation of 3.8, if $L_C$ is decidable, then $C$ is decidable. For example, if $X$ is a perfectly normal Hausdorff space that has a decidable lattice of cozero sets, then the $\ell$-group $C$ of continuous functions $X \longrightarrow \mathbb{R}$ is decidable.*

*Proof.* After attaching Gödel numbers to formulas of $\mathscr{L}^+$ and $\mathscr{L}^{lat}$ it is a routine matter to show that there is a recursive function $f : \mathbb{N} \longrightarrow \mathbb{N}$ with $\ulcorner\varphi_\Lambda\urcorner = f(\ulcorner\varphi\urcorner)$ for all $\mathscr{L}^+$-formulas $\varphi$. Using 3.8 for sentences we get the assertion. $\square$

**3.10. Corollary.** *Let $M \prec N$ be an elementary extension of o-minimal expansions of DOAGs (in the same signature) and let $\overline{N}$ be an o-minimal expansion of $N$. Let $X \subseteq M^n$ be $M$-definable. Let*

$$C = \{f : X \longrightarrow M \mid f \text{ is continuous and } M\text{-definable}\}$$
$$D = \{f : X_N \longrightarrow N \mid f \text{ is continuous and } \overline{N}\text{-definable}\}.$$

*(Here $X_N$ denotes the subset of $N^n$ defined in $N$ by a formula that defines $X$.) We consider $C$ and $D$ as $\ell$-groups. Assume $C$ and $D$ are constructible $\ell$-groups, see 2.2 for conditions when this is the case. Using distance functions we know that $L_C$ consists of the $M$-definable subsets of $X$ that are open in $X$ and $L_D$ consists of the $\overline{N}$-definable subsets of $X_N$ that are open in $X_N$.*

*Obviously the map $C \longrightarrow D$ that sends $f$ to $f_N$ is an embedding of $\ell$-groups and the map $L \longrightarrow L'$ that sends $U$ to $U_N$ is an embedding of bounded lattices.*

*If the map $L_C \longrightarrow L_D$ is elementary, then so is the map $C \longrightarrow D$. In particular, if $N = M$ and $X$ is of dimension 1, then the map $C \longrightarrow D$ is elementary.*

*Proof.* This is immediate from 3.8. Notice that in the case $\dim(X) = 1$ we have $L_D = L_C$ by o-minimality. $\square$

---

[9] The right hand side has to be understood as substitution of the variables $v_{s,t}$ by $\{s(\bar{f}) < t(\bar{f})\}$, which is the entry of $\Lambda_{d,\bar{x}}(\bar{f})$ at the index $(s,t) \in \text{Ht}_{\bar{x}}(d)^2$; then by 3.7, the right hand side of the equivalence indeed has a truth value.

4. APPLICATIONS

Before applying the results of the previous section we need to talk about some model theory of topological lattices. This will mostly be a recap of known results and some consequences of these, adapted to our needs. Decidability of some of these lattices will be imported from monadic second order logic.

4.1. **Some model theory of topological lattices**

(i) Let $M$ be a first order structure in a language $\mathscr{L}$. The **monadic second order structure of** $M$ is defined to be the following first order structure $\mathrm{MSO}(M)$: The universe is the powerset of $M$. Then $\mathrm{MSO}(M)$ is the expansion of the partially ordered set given by inclusion of subsets of $M$ by the structure that is induced by $M$ on the atoms of $(\mathrm{MSO}(M), \subseteq)$ via the bijection $M \longrightarrow \mathrm{Atoms}(\mathrm{MSO}(M))$, $a \mapsto \{a\}$. We generally identify the elements of $M$ with the atoms of $\mathrm{MSO}(M)$. Then every subset $S \subseteq M^n$, 0-definable in $M$ is also 0-definable in $\mathrm{MSO}(M)$.

The **weak monadic second order structure of** $M$ is defined to be the substructure $W(M)$ of $\mathrm{MSO}(M)$ induced on the finite subsets of $M$. Explicitly, $W(M)$ is the first order structure expanding the partially ordered set of finite subsets of $M$ by the 0-definable (in $M$) subsets of $M^n$.

(ii) Let $S2$ be the binary tree $2^{<\omega}$ together with the two successor functions $\sigma \mapsto \sigma\hat{\ }1$ and $\sigma \mapsto \sigma\hat{\ }0$. Then any expansion of $\mathrm{MSO}(S2)$ by naming finitely many elements from $W(S2)$ is decidable.

This is the main result in [Rabin1969], see [Rabin1969, Theorem 1.1] and [Rabin1969, Corollary 1.9].

(iii) Let $S \prec T$ be an elementary extension of densely ordered sets. Then $W(S) \prec W(T)$. Since densely ordered sets are $\omega$-categorical, this follows from [Tressl, 2.5].

(iv) If $T$ is a densely linearly ordered set, then any expansion of $W(T)$ by naming finitely many elements is decidable, see [Tressl, 3.6], [Laeuch1968].

We turn our attention to lattices. Let $T$ be a dense linear order and let $L(T)$ be the bounded sublattice of the powerset of $T$ generated by the sets $(-\infty, a]$ and $[a, +\infty)$ with $a \in T$.

(v) The partially ordered set $(L(T), \subseteq)$ is interpretable in the first order structure $W(T)$. After naming 2 parameters, $(L(T), \subseteq)$ and $W(T)$ are bi-interpretable (see [Hodges1993, section 5.4, item (c), p. 222] for the meaning of this). This is [Tressl, 3.2].

(vi) Let $M$ be an expansion of a DOAG and let $X, Y \subseteq M^n$ be definable. We write $L(X)$ for the lattice of closed and definable subsets of $X$. Suppose $X$ and $Y$ are closed subsets of $S := X \cup Y$. Let $Z = X \cap Y$. Then the lattice $L(S)$ is isomorphic to the fibre product

$$L(X) \times_{L(Z)} L(Y) = \{(A, B) \in L(X) \times L(Y) \mid A \cap Z = B \cap Z\}$$

of the lattices $L(X)$ and $L(Y)$ along the natural maps $L(X) \twoheadrightarrow L(Z)$, $L(Y) \twoheadrightarrow L(Z)$ given by the inclusions $Z \hookrightarrow X$, $Z \hookrightarrow Y$. Consequently, $L(S)$ is definable in the expansion of the lattice $L(X) \times L(Y)$ by a unary predicate naming the diagonal $\Delta = \{(C, C) \mid C \in L(Z)\}$ of $L(Z)$, using $(Z, Z)$ as a parameter:

$$A \cap Z = B \cap Z \iff (A, B) \wedge (Z, Z) \in \Delta.$$

If $Z$ is finite, then $\Delta$ is certainly parametrically definable in $L(X) \times L(Y)$ and so $L(S)$ is parametrically definable in the lattice $L(X) \times L(Y)$.

(vii) Let $M$ be an o-minimal structure and let $X \subseteq M^n$ be definable and of dimension 1.
   (a) The lattice $L(X)$ of closed definable subsets of $X$ is definable with parameters in the lattice $L(M)$.[10]
   (b) If $M \prec N$, then $L(X) \prec L(X_N)$.
   (c) Any expansion of the lattice $L(X)$ by finitely many parameters is decidable.

*Proof.* (a) This follows from (vi), since by o-minimality each infinite definably connected component $X_0$ of $X$ is obtained from $M$ by glueing (in the sense of (vi)) intervals (better: definable homeomorphic copies of intervals) with endpoints in $M$ along finite sets. The lattice $L(X)$ is then the product of the $L(X_0)$.

(b) The reasoning in (a) is natural in $M$ and thus $L(X) \prec L(X_N)$ follows from $L(M) \prec L(N)$, which itself is a consequence of (v) and (iii).

(c) follows from (a) by (iv) and (v). □

It should be said at this point that by strengthening techniques from [Grzego1951] on undecidable Heyting algebras one can show that the lattice of all closed 1-dimensional semilinear subsets of $K^2$ interprets $(\mathbb{N}, +, \cdot)$ and is thus undecidable, see [Tressl, 6.2]

4.2. **Application I**

Let $M$ be an o-minimal expansion of a DOAG such that for all $x, a, y, b \in M$ with $x \neq y$ there is a continuous $M$-definable function $M \longrightarrow M$ that maps $x$ to $a$ and $y$ to $b$. Then for every 1-dimensional definable subset $X$ of $M^n$, the $\ell$-group $C$ of continuous $M$-definable functions $X \longrightarrow M$ is decidable.

For example, the $\mathbb{Q}$-vector lattice $V$ of all semilinear continuous functions $S^1 \longrightarrow \mathbb{Q}$ is decidable, where $S^1$ denotes the set of all $(x, y) \in \mathbb{Q}^2$ with $\max\{|x|, |y|\} = 1$. It is well known that $V$ is the free vector lattice in two generators, cf. [Beynon1974].

*Proof.* From the assumption on $M$ it follows easily that property (C3) holds for 2-element subsets as required in 2.2(ii) and therefore $C$ is a constructible $\ell$-group. By 4.1(vii), the lattice $L(X)$ is decidable and hence also its inverse $L_C$ is decidable. By 3.9, $C$ is decidable. □

4.3. *Remarks.*

(i) Suppose $K$ is an ordered field and $M$ is the ordered $K$-vector space $K$. Let $X \subseteq M^n$ be semi-linear and of dimension 1. By 4.2, $C$ is decidable and this is built on the decidability of the weak monadic second order theory of $(\mathbb{Q}, \leq)$. One might wonder if this is necessary; in fact it is:

It is not difficult to show that the lattice of (co-)zero sets of $X$ is interpretable in $C$, provided the set $P$ of all $f \in C$ that have no zeroes is definable in $C$. For example for $X = K$ this is an easy exercise, for $X = S^1$ this is done in [GlMaPo2005, Lemma 4.4(3)] (in fact the key point in section 4 of [GlMaPo2005] is to show that $P$ is definable even if $X = S^n$). Once the lattice of zero sets of $X$ is interpreted one can also interpret the lattice of closed definable subsets of $K$ and then by 4.1(v) we see that $C$ interprets the weak monadic second order theory of $(\mathbb{Q}, \leq)$. Hence any proof of the decidability of $C$ will also have to address the decidability of the weak monadic second order theory of $(\mathbb{Q}, \leq)$.
 (ii) The decidability of the free vector lattice in two generators is claimed in [GlMaPo2005, Theorem 5.1] with an incorrect proof, see [GlMaPo2016]. Hence 4.2 reinstates [GlMaPo2005, Theorem 5.1]. On the other hand the decidability of the free abelian $\ell$-group in 2 generators remains open. One encounters two problems when trying to apply the method from sections 2 and 3: Firstly, 2.3 – as it stands – is not available without division by integers (Notice that the proof of 2.4 only uses divisibility of the $\ell$-group $C$ in that it invokes 2.3.) The second problem comes with 3.4: Linear equations in $x_1, \ldots, x_n, y$ are not explicit in $y$ unless we allow division by the coefficient of $y$.

---

[10]If $K$ is an ordered field and $M$ expands the ordered $K$-vector space $K$, then $L(M)$ is in fact parametrically bi-interpretable with $L(X)$. We leave this argument to the reader.

### 4.4. Application II
The $\ell$-group $C$ of continuous functions $\mathbb{R} \longrightarrow \mathbb{R}$ is decidable.

*Proof.* By 3.9 it suffices to show that the lattice $L_C$ is decidable. Now $L_C$ is the lattice of open subsets of $\mathbb{R}$ and this lattice is decidable as a consequence of 4.1(ii), see [Rabin1969, Theorem 2.9]. □

We now turn to model theoretic applications, using results from [Astier2013]:

### 4.5. Theorem.
 (i) If $M$ is an o-minimal expansion of a real closed field and $X \subseteq M^n$ is semi-linear and bounded, then the lattice of semi-linear subsets of $X$ that are open in $X$ is an elementary substructure of the lattice of $M$-definable subsets of $X$ that are open in $X$. This follows from [Astier2013, Corollary 3.7 item 2.].
 (ii) If $M \prec N$ are o-minimal expansions of real closed fields and $X \subseteq M^n$ is $M$-definable, then the lattice of $M$-definable subsets of $X$ that are open in $X$ is an elementary substructure of the lattice of $N$-definable subsets of $X_N$ that are open in $X_N$. This follows from [Astier2013, Corollary 2.16 item 2.].

### 4.6. Application III
Let $M \prec N$ be an o-minimal expansions of DOAGs and let $\overline{N}$ be an o-minimal expansion of $N$. Let $X \subseteq M^n$ be definable, let $C$ be the $\ell$-group of continuous $M$-definable functions $X \longrightarrow M$ and let $D$ be the $\ell$-group of continuous $\overline{N}$-definable functions $X_N \longrightarrow N$.
 (i) If $X$ is of dimension 1 and $M$ expands the ordered $K$-vector spaces of an ordered field $K$ and $\overline{N}$ expands the ordered $L$-vector spaces of an ordered field $L \supseteq K$, then the natural map $C \longrightarrow D$ is elementary.

Now assume $M$ and $N$ expand fields and Let $D_0 \subseteq D$ be the $\ell$-group of continuous $N$-definable functions $X_N \longrightarrow N$.
 (ii) The natural map $C \longrightarrow D_0$ is elementary.
 (iii) $C$ is elementary equivalent to $D$.
 (iv) If $X$ is bounded and semi-linear, then the inclusion $C_{s.l.}(X) \hookrightarrow C$ and the natural map $C_{s.l.}(X) \longrightarrow C_{s.l.}(X_N)$ are elementary.
 (v) Every convex lattice ordered subgroup of $C$ containing a nonzero constant function is an elementary substructure in $C$.

*Proof.* (i) By 2.2(ii), $C$ and $D$ are constructible. Hence by 3.10 it suffices to show that $L_C \prec L_D$. This follows from 4.1(vii), since the lattice of closed and $N$-definable subsets of $X_N$ is the lattice of closed and $\overline{N}$-definable subsets of $X_N$.

(ii), (iii) and (iv). Now assume $M$ and $N$ expand fields. By 2.2(iii),(iv) all the $\ell$-groups in (ii),(iii) and (iv) are constructible. Hence by 3.10 again it suffices to show that the corresponding maps of lattices of cozero sets in (ii) and (iv) are elementary. Item (iii) is a consequence of (ii) and (iv).

As $M$ expands a field, $X$ is $M$-definably homeomorphic to a bounded semi-linear set by the triangulation theorem in o-minimal expansions of fields, cf. [vdDries1998, Chapter 8, (1.7), p. 122]. Hence the corresponding $\ell$-group of continuous definable functions are naturally isomorphic and we may assume that $X$ is a bounded semi-linear set in (ii), too. Using 4.5(i), we see that the maps $C_{s.l.}(X) \hookrightarrow C$ and $C_{s.l.}(X_N) \hookrightarrow D_0$ are elementary. Using 4.5(ii), we see that the map $C \longrightarrow D_0$ is elementary. But then the map $C_{s.l.}(X) \longrightarrow C_{s.l.}(X_N)$ is elementary, too. This shows (ii) and (iv).

(v) Let $C_0$ be a convex lattice ordered subgroup of $C$ containing a nonzero constant function. It is easy to see that $C_0$ is again a constructible $\ell$-group with the same zero sets as $C$. Hence 3.10 applies.    □

4.7. *Example.* Neither the assumption that $M$ expands a field for 4.6(v), nor the assumption that $X$ is bounded in 4.6(iv) can be dropped. If $K$ is any ordered field then the $\ell$-group $G$ of semi-linear continuous functions $K^2 \longrightarrow K$ fails to have (C3). In fact, $G$ has a convex lattice ordered subgroup with the same zero sets as $G$, but which is not an elementary substructure of $G$:

By [AscTha2015, Example 3.4] the semi-linear functions $y$ defined on $\{x \geq 1\} \subseteq K^2$ and $-y$ defined on $\{x \leq 0\} \subseteq K^2$ have no continuous semi-linear extension $K^2 \longrightarrow K$. (Notice that the two sets are actually disjoint.)

A slight variation of this example has an interesting application in our context. Let $H$ be the open left half plane $\{(x,y) \in K^2 \mid x < 0\}$ and let $f_0 \in G$ be the function $(0 \vee -x) \wedge 1$. Hence $f_0 \geq 0$ is bounded and $H$ is the cozero set of $f_0$. If $f \in G$, then $H \subseteq \{f = 0\}$ if and only if $|f| \wedge f_0 = 0$.

Further, let $S$ be the open vertical strip $\{(x,y) \in K^2 \mid 1 < x < 2\}$ and let $g_0 = ((x-1) \vee 0) \wedge ((2-x) \vee 0) \in G$. Then $g_0 \geq 0$ is bounded and $S$ is the cozero set of $g_0$. If $f \in G$, then $S \subseteq \{f = 0\}$ if and only if $|f| \wedge g_0 = 0$.

Now let $G_0 = \{f \in G \mid f$ is bounded in $K$ on $S\}$. Obviously $G_0$ is a convex lattice ordered subgroup of $G$ having the same zero sets as $G$. Let $\varphi(f,g)$ be the formula (with parameters $f_0, g_0 \in G_0$)

$$||f - |g|| \wedge g_0 = 0 \ \& \ |f| \wedge |f_0| = 0,$$

expressing that

$$S \subseteq \{f = |g|\} \ \& \ H \subseteq \{f = 0\}.$$

Then $G_0 \models \forall g \exists f\ \varphi(f,g)$ and $G_0$ is defined in $G$ by $\exists f\ \varphi(f,g)$. In particular $G_0$ is not an elementary substructure of $G$.

*Proof.* If $g \in G_0$, say $|g| \leq N \in K$ on $S$, then $f = (0 \vee N \cdot x) \wedge |g| \in G_0$ satisfies $\varphi(f,g)$. Conversely, if $g \notin G_0$, then $|g|$ is unbounded on $S$ and the argument in

[AscTha2015, Example 3.4] shows that there is no $f \in G$ that vanishes on $H$ and extends $|g|$ on $S$.[11] □

Thus $G_0 \not\prec G$, which shows that 4.6(v) fails without the assumption on $M$ expanding a field. Now assume $K$ is a real closed field and $M$ is this field. Then in $C$, property (C3) holds, and therefore $C \models \forall g \exists f \ \varphi(f,g)$. Hence $G = C_{s.l.}(K^2)$ is not an elementary substructure of $C$, which shows that the boundedness assumption in 4.6(iv) cannot be dropped.

### 4.8. Application IV

Let $M$ be an o-minimal expansion of a real closed field and let $X \subseteq M^n$ be $M$-definable. The following $\ell$-groups are elementary equivalent:

(i) The $\ell$-group $C$ of continuous functions $X \longrightarrow M$ for addition.
(ii) The $\ell$-group $D$ of continuous functions $X \longrightarrow (0,\infty)_M$ for multiplication.

(This statement is only of relevance if $M$ does **not** expand a model of real exponentiation.)

*Proof.* Notice that $C$ and $D$ have the same cozero sets: $L_D$ consist of all open sets of the form $\{f \neq 1\}$ where $f \in D$. But these are obviously just the open $M$-definable subsets of $X$. Hence by 3.8, $C$ and $D$ are elementarily equivalent (3.8 is applicable by 2.2(iv),(v)). □

---

[11] or one can deduce it directly from the semilinear Pierce-Birkhoff representation of $f \in G$, which says that each $f$ is a finite supremum of finite infima of affine functions $K^2 \longrightarrow K$; hence the partial derivative $\frac{\partial}{\partial x} f$ (where defined) is bounded in $K$ and the continuous and piecewise differentiable map $f$ cannot vanish on the $y$-axis, if it is unbounded in $K$ on $S$.

The University of Manchester, School of Mathematics, Oxford Road, Manchester M13 9PL, UK

Homepage: http://personalpages.manchester.ac.uk/staff/Marcus.Tressl/index.php

*E-mail address*: marcus.tressl@manchester.ac.uk